\documentclass[12pt]{article}
\usepackage[final]{epsfig}	
\usepackage{graphics}
\usepackage{amsmath}
\usepackage{amsfonts}
\usepackage{latexsym}
\usepackage{amssymb}

\newtheorem{lemma}{Lemma}[section]

\newtheorem{remark}[lemma]{Remark}

\newtheorem{theorem}{Theorem}

\newtheorem{corollary}[theorem]{Corollary}

\begin{document}
\newcommand{\eps}{{\varepsilon}}
\newcommand{\proofend}{$\Box$\bigskip}
\newcommand{\C}{{\mathbf C}}
\newcommand{\Q}{{\mathbf Q}}
\newcommand{\R}{{\mathbf R}}
\newcommand{\Z}{{\mathbf Z}}
\newcommand{\RP}{{\mathbf {RP}}}

\title {The Poncelet grid and the billiard in an ellipse}
\author{Mark Levi  and Serge Tabachnikov\\
{\it Department of Mathematics, Penn State University}\\
{\it University Park, PA 16802, USA}}
\date{}
\maketitle

\section{The closure theorem and the Poncelet grid} \label{closure}

\begin{figure}[ht] 
\centerline{\epsfbox{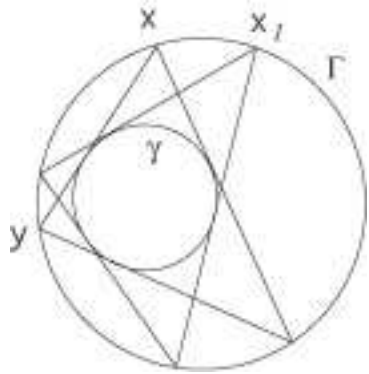}}
\caption{Poncelet polygons}
\label{clos}
\end{figure}

The Poncelet closure theorem (or Poncelet porism) is a classical result of projective geometry. Given two nested ellipses, $\gamma$ and $\Gamma$, one plays the following game: choose a point $x$ on $\Gamma$, draw a tangent line to $\gamma$ until it intersects $\Gamma$ at point $y$, repeat the construction, starting with $y$, and so on. One obtains a polygonal line, inscribed into $\Gamma$ and circumscribed about $\gamma$. Suppose that this process is periodic: the $n$-th point coincides with the initial one.  Now start at a different point, say, $x_1$. The Poncelet closure theorem states that   the polygonal line again closes up after $n$ steps, see figure \ref{clos}.  We will call these closed inscribed-circumscribed lines Poncelet polygons.

Although the Poncelet theorem is almost 200 years old, it continues to attract interest: see \cite{B-B, Ber, B-K-O-R, G-H, Mir1, Mir2, Tab93, Tab95} for a sample of references. 

\begin{figure}[ht]
\centerline{\epsfbox{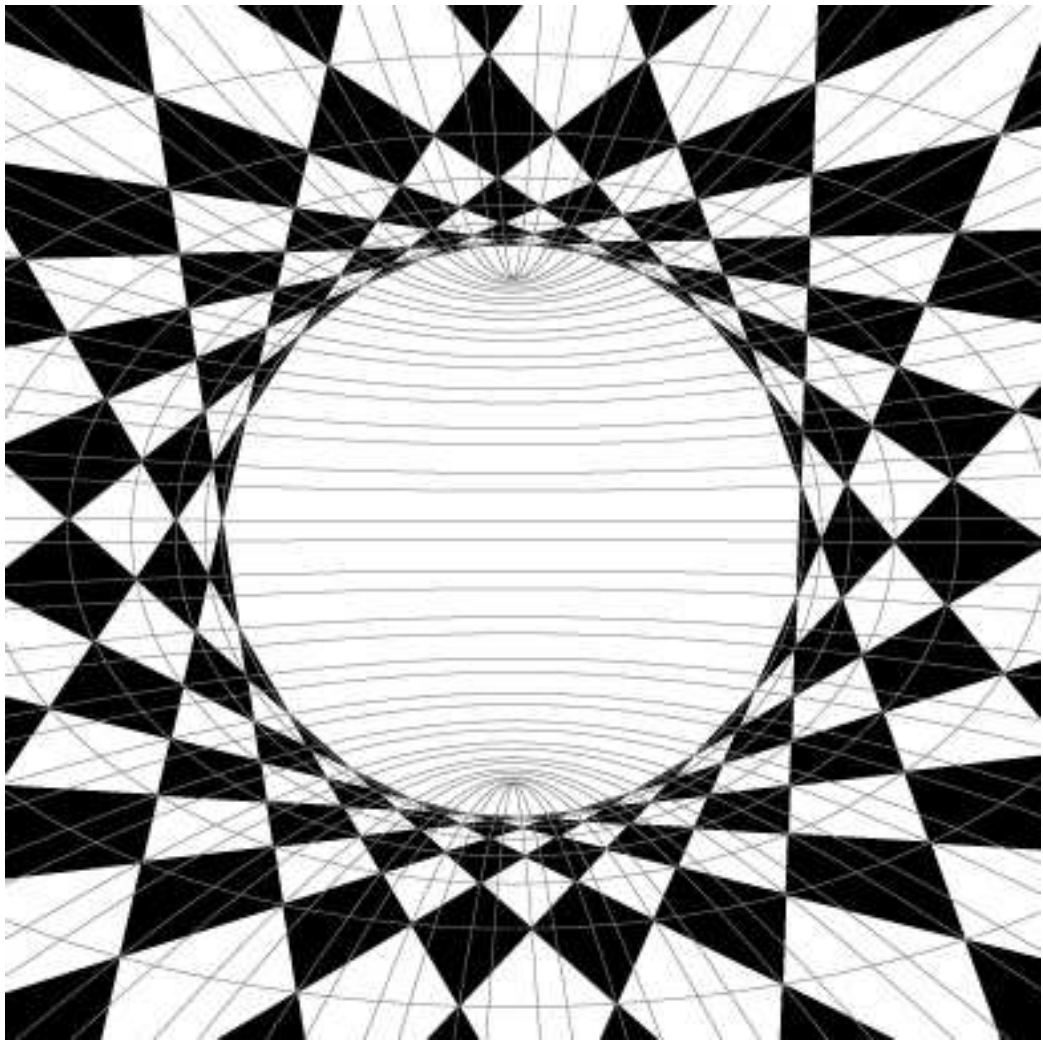}}
\caption{Poncelet grid}
\label{gr}
\end{figure}

Recently R. Schwartz \cite{Sch} discovered the following property of Poncelet polygons. Extending the sides of a Poncelet $n$-gon, one obtains a set of points  called the Poncelet grid, see figure \ref{gr} borrowed from \cite{Sch}. The points of the Poncelet grid can be viewed as lying on a family of nested closed curves, and also on a family of disjoint curves having radial directions. 

More precisely, let $\ell_1,\dots,\ell_n$ be the lines containing the sides of the polygon, enumerated in such a way that their tangency points to $\gamma$ are in cyclic order. 
The Poncelet grid consists of $n(n+1)/2$ points $\ell_i \cap \ell_j$. The indices are understood cyclically and, by convention, $\ell_j \cap \ell_j$ is the tangency point of $\ell_j$ with $\gamma$.
Define the sets: 
\begin{equation} \label{sets}
P_k = \cup_{i-j=k} \ell_i \cap \ell_j,\quad Q_k = \cup_{i+j=k} \ell_i \cap \ell_j.
\end{equation}
The cases of odd and even $n$ differ somewhat and, as in \cite{Sch}, we assume that $n$ is odd. There are $(n+1)/2$ sets $P_k$, each containing $n$ points, and $n$ sets $Q_k$, each containing $(n+1)/2$ points. 

The Schwartz theorem states: 

\begin{theorem} \label{RSCH}
The sets $P_k$ lie on nested ellipses, and the sets $Q_k$ on disjoint hyperbolas; the complexified versions of these ellipses and hyperbolas have four common complex tangent lines. Furthermore,  all the sets $P$s are projectively equivalent to each other, and all the sets $Q$s are projectively equivalent to each other. 
\end{theorem}

The proof in \cite{Sch} consists in a study of  properties of the underlining elliptic curve; we will give a different, more elementary, proof and deduce this theorem from properties of  billiards in ellipses.

\section{Mathematical billiards: general facts} \label{bill}
In this section we recall (with proofs) necessary facts about billiards, see \cite{Tab95, Tab05} for detailed surveys. 

The billiard system describes the motion of a free point inside a plane domain: the point moves with a constant speed along a straight line until it hits the boundary, where it reflects according to the familiar law of geometrical optics ``the angle of incidence equals the angle of reflection". 

We assume that the billiard table is a convex domain with a smooth boundary curve $\Gamma$. The billiard ball map acts on oriented lines that intersect the billiard table, sending the incoming billiard trajectory to the outgoing one. Let $x,y,z$ be points on $\Gamma$ such that the line $xy$ reflects to the line $yz$. The equal angles condition has a variational meaning.

\begin{lemma} \label{var}
The angles made by lines $xy$ and $yz$ with $\Gamma$ are equal if and only if $y$ is a critical point of the function $|xy|+|yz|$ (considering $x$ and $z$ as fixed).
\end{lemma} 

\noindent{\bf Proof}. Assume first that $y$ is a free point, not confined to $\Gamma$. The gradient of the function $|xy|$ is the unit vector from $x$ to $y$, and the gradient of $|yz|$ is the unit vector from $z$ to $y$. By the Lagrange multipliers principle,  $y\in \Gamma$ is a critical point of the function $|xy|+|yz|$ if and only if the sum of the two gradients is orthogonal to $\Gamma$, and this is equivalent to the fact that $xy$ and $yz$ make equal angles with $\Gamma$.
\proofend

An important consequence is that the billiard ball map is area preserving. We continue to identify oriented lines intersecting $\Gamma$ with pairs of points  $(x,y)$.

\begin{theorem} \label{area}
The area element 
$$
\omega=\frac{\partial^2 |xy|}{\partial x \partial y}\ dx \wedge dy
$$
is invariant under the billiard ball map.
\end{theorem}

\noindent{\bf Proof}. According to Lemma \ref{var}, 
$$
\frac{\partial |xy|}{ \partial y} + \frac{\partial |yz|}{\partial y}=0.
$$
Take the differential:
$$
\frac{\partial^2 |xy|}{\partial x \partial y}\ dx + \frac{\partial^2 |xy|}{\partial y^2}\ dy + \frac{\partial^2 |yz|}{\partial y^2}\ dy + \frac{\partial^2 |yz|}{\partial y \partial z}\ dz=0,
$$
and wedge multiply by $dy$ to obtain
$$
\frac{\partial^2 |xy|}{\partial x \partial y}\ dx \wedge dy = \frac{\partial^2 |yz|}{\partial y \partial z}\ dy \wedge dz,
$$
as claimed.
\proofend

\begin{remark} {\rm One can express the invariant area form $\omega$ in more convenient coordinates. Although it is of no importance to us, let us mention two such formulas. First, one can characterize an oriented line by a point $x \in \Gamma$ and the angle $\alpha$ made with $\Gamma$ at $x$. Assuming that $x$ is an arc length parameter on $\Gamma$, one has: $\omega=\sin \alpha\ d\alpha \wedge dx$. We will provide a geometrical explanation of Theorem \ref{area} in an Appendix. Secondly, an oriented line is characterized by its direction $\varphi$ and its signed distance $p$ from an origin. Then  $\omega=dp \wedge d\varphi$.}
\end{remark}

Another necessary fact about billiards  concerns caustics. A  caustic  is a curve inside a  billiard table such that if a segment of a billiard trajectory is tangent to this curve, then so is each reflected segment. We assume that caustics are smooth and convex.

Let  $\Gamma$ be the boundary of a billiard table and $\gamma$ a caustic. Suppose that one erases the table, and only the caustic  remains. Can one recover $\Gamma$ from $\gamma$? The answer is  given by the following {\em string
construction}: wrap a closed non-stretchable string around $\gamma$, pull it tight at a
point and move this point around $\gamma$ to obtain a curve $\Gamma$.

\begin{theorem} \label{string}
The billiard inside $\Gamma$ has $\gamma$ as its caustic.
\end{theorem}

\begin{figure}[ht]
\centerline{\epsfbox{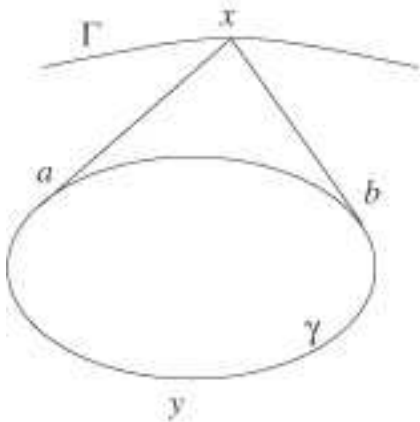}}
\caption{String construction}
\label{stringconstr}
\end{figure}

\noindent{\bf Proof}.  Choose a reference point $y\in \gamma$. For a point $x\in \Gamma$, let $f(x)$ and $g(x)$ be the distances from $x$ to $y$ by going around $\gamma$ on the right and on the left, respectively. Then $\Gamma$ is a level curve of the function $f+g$. We want to prove that the angles made by the segments $ax$ and $bx$ with $\Gamma$ are equal; see figure~\ref{stringconstr}.

We claim that  the gradient of $f$ at $x$  is the unit vector in the direction $ax$.
Indeed,  the free end $x$ of the contracting string $yax$  will move directly toward point $a$ with unit speed. 
It follows  that $\nabla(f+g)$ bisects the angle $axb$.
Therefore $ax$ and $bx$ make equal angles with $\Gamma$.
\proofend

Note that the string construction provides a one-parameter family of billiard tables: the parameter is the length of the string. Note also that, by the same reasoning,  the level curve of the function $f-g$ is orthogonal to $\Gamma$.

\section{The billiard in an ellipse: integrability and its consequences} \label{ell}
Optical properties of conics were already known to the Ancient Greeks. In this section we review billiards in ellipses and describe some consequences of their complete integrability.

\begin{figure}[ht]
\centerline{\epsfbox{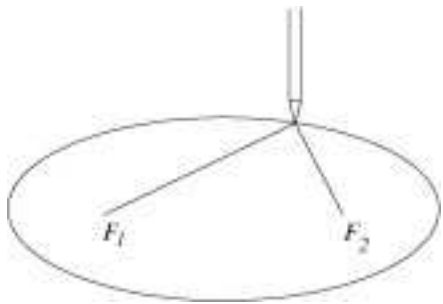}}
\caption{Gardener's construction of an ellipse}
\label{gard}
\end{figure}

First of all, recall the geometric definition of an ellipse: it is the locus of points whose sum of distances to two given points, $F_1$ and $F_2$, is fixed; these two points are called the foci. An ellipse can be constructed using a string whose ends are fixed at the foci, see figure \ref{gard}.  A 
hyperbola is defined similarly with the sum of distances replaced by the absolute
value of their difference. Taking the segment $F_1 F_2$ as $\gamma$ in Theorem \ref{string}, it follows that a ray passing through one focus reflects to a ray passing through the other focus.

The construction of an ellipse with given foci has a parameter, the length of the
string. The family of conics with fixed foci is called confocal. The equation
of a confocal family, including ellipses and hyperbolas, is
\begin{equation} \label{conffam}
\frac{x_1^2}{a_1^2+\lambda}+\frac{x_2^2}{a_2^2+\lambda}=1
\end{equation}
where $\lambda$ is a parameter.

Fix $F_1$ and $F_2$. Given a generic point  in the plane, there exist a unique
ellipse and a unique hyperbola with foci $F_1, F_2$ passing through the point. The ellipse and the hyperbola are orthogonal to each other: this follows from the fact that the sum of two unit vectors is perpendicular to its difference; cf. proofs of Lemma \ref{var} and Theorem \ref{string}. 

The next theorem says that the billiard ball map  in an ellipse is  integrable, that is, possesses  an invariant quantity. 

\begin{theorem} \label{integr} 
A billiard trajectory inside an ellipse forever remains tangent to a fixed confocal conic. More precisely, if a segment of a billiard trajectory does not intersect the  segment $F_1 F_2$, then all the segments of this trajectory do not intersect $F_1 F_2$ and are all tangent to the same ellipse with  foci $F_1$ and $F_2$; and if a segment of  a trajectory intersects $F_1 F_2$, then all the segments of this trajectory intersect $F_1 F_2$ and  are all tangent to the same hyperbola with  foci $F_1$ and $F_2$. 
\end{theorem}

Thus the billiard inside an ellipse has a 1-parameter family of caustics consisting of 
confocal ellipses. We give an elementary geometry proof of Theorem \ref{integr} in  an Appendix.

Theorems \ref{string} and \ref{integr} imply the following Graves theorem: {\it wrapping a closed non-stretchable string around an ellipse produces a confocal ellipse}, see \cite{Ber, Po}.

\begin{figure}[ht]
\centerline{\epsfbox{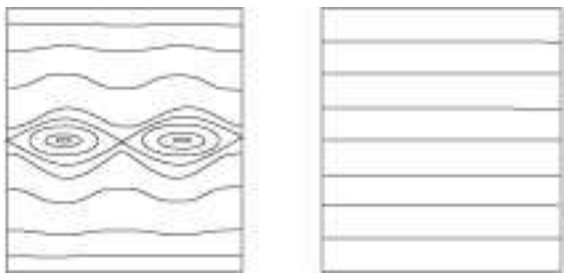}}
\caption{Phase space of the billiard ball map in an ellipse and in a circle}
\label{phase}
\end{figure}

The space of oriented lines intersecting an ellipse is, topologically, a cylinder. This cylinder is foliated by invariant curves of the billiard ball map, see figure \ref{phase} on the left. 
Each curve represents the family of rays tangent to a fixed confocal conic. The $\infty$-shaped curve corresponds to the family of rays through the foci. The two singular points of this curve represent the major axis with two opposite orientations, a 2-periodic
billiard trajectory. Another 2-periodic trajectory is the minor axis represented by
two centers of the regions inside the $\infty$-shaped curve. For comparison, we also give a phase portrait of the billiard ball map in a circle, see figure \ref{phase} on the right.

The integrability of the billiard ball map makes it possible to choose a cyclic coordinate on each invariant curve, say, $x$ modd 1, such that the map is given by a shift $x\mapsto x+c$; the value of the constant $c$ depends on the invariant curve. Here is a description of this construction.

Choose a function $f$ on the cylinder whose level curves are the invariant curves of the billiard ball map. Let $\gamma$ be a curve $f=c$. Consider the curve $\gamma_{\varepsilon}$ given by $f=c+\varepsilon$. For an  interval $I \subset \gamma$, consider the area $\omega(I,\varepsilon)$ between $\gamma$ and $\gamma_{\varepsilon}$
over $I$. Define the ``length" of $I$ as 
$$
\lim_{\varepsilon\to 0} \frac{\omega(I,\varepsilon)}{\varepsilon}.
$$
Choosing a different function $f$, one multiplies the length of every segment by the
same factor. Choose a coordinate $x$ so that the length element is $dx$; this
coordinate is well defined up to an affine transformation.  Normalizing $x$ so that the total length is 1 determines $x$ up to a shift $x\mapsto x+const$.

The billiard ball map preserves the area element $\omega$ and the invariant curves. Therefore it preserves the length element on the invariant curves, that is, is given by the formula $x\mapsto x+c$.

Let us summarize. Consider an ellipse $\Gamma$ and a confocal ellipse $\gamma$, a caustic for the billiard in $\Gamma$. The billiard ball map is a self-map of $\gamma$ (it sends point $a$ to $b$ in figure \ref{stringconstr}). We have introduced a parameter $x$ on $\gamma$ such that the billiard ball map is a shift $x\mapsto x+c$; the value of $c$ depends on $\Gamma$, but the parameter $x$ depends  on $\gamma$ only.

Let $\Gamma'$ be another confocal ellipse containing $\gamma$. Then $\Gamma$ and $\Gamma'$ share the caustics, in particular, $\gamma$. It follows that the billiard ball map associated with $\Gamma'$ is also a shift  in the parameter $x$.

\begin{figure}[ht]
\centerline{\epsfbox{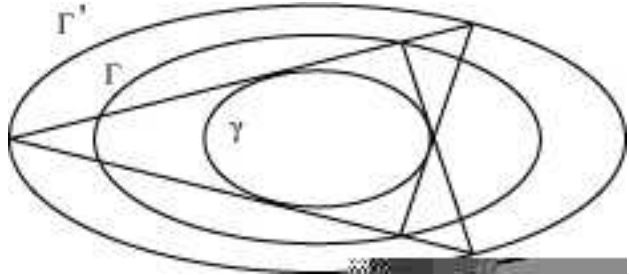}}
\caption{Commuting billiard ball maps}
\label{commute}
\end{figure}

\begin{corollary} \label{comm}
The billiard ball maps associated with $\Gamma$ and $\Gamma'$ commute, see figure \ref{commute}.
\end{corollary} 

\noindent{\bf Proof}. The shifts $x\mapsto x+c$ and $x\mapsto x+c'$commute.
\proofend

Another consequence is a Poncelet-style closure theorem.

\begin{corollary} \label{closcor}
Assume that a billiard trajectory in an ellipse $\Gamma$, tangent to a confocal ellipse $\gamma$, is $n$-periodic. Then every billiard trajectory in $\Gamma$, tangent to $\gamma$, is $n$-periodic. 
\end{corollary} 

\noindent{\bf Proof}. In the appropriate coordinate on $\gamma$, the billiard ball map
is $x\mapsto x+c$. A point is $n$-periodic if and only if $nc$ is an integer. This condition does not depend on $x$, and the result follows.
\proofend

One can further generalize. Let $\Gamma_1,\Gamma_2,\dots,\Gamma_k$ be confocal ellipses and $\gamma$ another confocal ellipse inside them all. One may modify the formulation of Corollary \ref{closcor} replacing a single billiard ball map by the composition of the billiard ball maps associated with $\Gamma_1,\Gamma_2,\dots,\Gamma_k$: the conclusion of the closure
theorem will hold without change.

Finally, being only a particular case of the Poncelet porism, Corollary \ref{closcor} implies its general  version. This is because a generic pair of nested ellipses is projectively equivalent to a pair of confocal ones (this proof of the Poncelet porism is mentioned in \cite{Ve}). 

More precisely, consider the complexified situation. Two conics have four common  tangent lines, and one has a 1-parameter family of conics sharing these four tangents. 

\begin{lemma} \label{confl}
A confocal family of conics consists of the conics, tangent to four fixed  lines.
\end{lemma}

\noindent{\bf Proof}. A curve, projectively dual to a conic, is a conic. The 1-parameter family of conics, dual to the confocal family (\ref{conffam}), is given by the equation 
$$
(a_1^2+\lambda) x_1^2+(a_2^2+\lambda) x_2^2=1.
$$
This is an equation of a pencil, a 1-parameter family of conics that pass through four fixed points; these are the  intersections of the two conics, $a_1^2 x_1^2+a_2^2 x_2^2=1$ and $ x_1^2+ x_2^2=1$. Projective duality interchanges points and tangent lines; applied again, it  yields a 1-parameter family of conics sharing four  tangent lines. 
\proofend

Since projective transformations act transitively of quadruples of lines in general position, 
a generic pair of conics is projectively equivalent to a pair of confocal ones.

\section{Back to the Poncelet grid} \label{back}
 Let $\gamma$ and $\Gamma$ be a pair of nested ellipses and $P$ a Poncelet $n$-gon circumscribing $\gamma$ and inscribed into $\Gamma$. Applying a projective transformation, we assume that $\gamma$ and $\Gamma$ are confocal.
 
 Let $x$ be the parameter on $\gamma$ introduced in Section \ref{ell}. Choosing the origin appropriately,   the consecutive tangency points of the sides of $P$ with $\gamma$ have coordinates 
 $$
 0,\ \frac{1}{n},\ \frac{2}{n}\ \dots,\ \frac{n-1}{n}.
 $$
 The set $P_k$ in (\ref{sets}) lies on the locus of intersections of the tangent lines to $\gamma$ at points $\gamma(x)$ and $\gamma(x+k/n)$ where $x$ varies from 0 to 1. This locus is a confocal ellipse for which  the billiard trajectories, tangent to $\gamma$, close up after $n$ reflections and $k$ turns around $\gamma$ (periodic trajectories with rotation number $k/n$). Thus $P_k$ lies on a confocal ellipse to $\gamma$.
 
Likewise, the set $Q_k$ in (\ref{sets}) lies on the locus of intersections of the tangent lines to $\gamma$ at points $\gamma(x)$ and $\gamma(k/n-x)$. We want to show that this locus is a confocal hyperbola. To this end we need the next result, which is an (apparently new) addition to  Theorem \ref{string}, the string construction. 

\begin{figure}[ht]
\centerline{\epsfbox{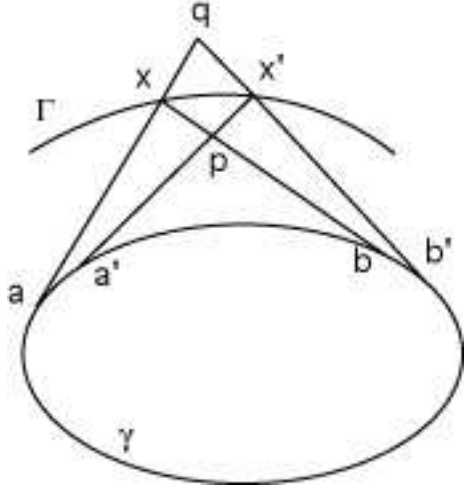}}
\caption{Addition to the string construction}
\label{strref}
\end{figure}

\begin{theorem} \label{addstr}
Apply the string construction to an oval $\gamma$ and let $x,x'$ be two infinitesimally close points on the curve $\Gamma$, see figure \ref{strref}. Then the line $pq$ is orthogonal to $\Gamma$.
\end{theorem}

\noindent{\bf Proof}. We will give two proofs, geometrical and analytical.
 
We claim that the infinitesimal quadrilateral $xqx'p$ is a rhombus. Indeed, by Theorem \ref{string}, the arc $xx'$ bisects the angles $pxq$ and $px'q$. Let $\varepsilon$ be the distance between $x$ and $x'$. Then the angles $pxq$ and $px'q$ are $\varepsilon$--close to each other. Now dilate with factor $1/\varepsilon$. The angles do not change, the arc $xx'$ becomes straight (to order $\varepsilon$), and, in the limit $\varepsilon \to 0$, we obtain a rhombus.

Analytically, let us see how fast the points $a$ and $b$ move as one moves the
point $x$ (not necessarily confined to $\Gamma$). Let the speeds of these points along $\gamma$ be $v_1$ and $v_2$; let the tangent segments $ax$ and $bx$ have lengths $l_1$ and $l_2$; let the angular velocity of the lines $ax$ and $bx$ be $\omega_1$ and $\omega_2$; and let $k_1$ and $k_2$ be the curvatures of $\gamma$ at points $a$ and $b$. Denote the velocity vector of point $x$ by $w$.

Then $k_1=\omega_1/v_1$ and $\omega_1=w_1/l_1$ where $w_1$ is the component of $w$ perpendicular to $ax$. Likewise, for the variables with index 2. It follows that 
$$
\frac{v_2}{v_1} = \frac{l_1k_1}{l_2k_2}\cdot \frac{w_2}{w_1}.
$$
 Consider two choices of $w$: tangent to
$\Gamma$ and perpendicular to it. Because of the equal angles property, Theorem \ref{string}, in the first case we have $w_1=w_2$, and in the second case, $w_1=-w_2$. Thus the ratio $v_2/v_1$ in both cases will have the same value and opposite
signs. This is equivalent to orthogonality of $xx'$ and $pq$.
 \proofend 
  
Now  we can describe the locus of intersections of the tangent lines to $\gamma$ at points $\gamma(x)$ and $\gamma(c-x)$. Indeed, by Theorem \ref{addstr}, this locus is a curve, orthogonal to the family of confocal ellipses, that is, a confocal hyperbola. It follows that the set $Q_k$ lies on a confocal hyperbola.
  
\section{Elliptic coordinates and linear equivalence of the sets $P$s and of the sets $Q$s}  \label{equiv}

It remains  to show that the sets $P_k$ are projectively (actually, linearly) equivalent for all values of $k$, and likewise for the sets $Q_k$.

Given an ellipse $\gamma$, let $x$ be the parameter on it described in Section \ref{ell}. 
Note that the map $x\mapsto x+1/2$ is central symmetry of the ellipse; in particular,
the tangent lines at points $\gamma(x)$ and $\gamma(x+1/2)$ are parallel.

For a point $P$ outside of $\gamma$, draw tangent segments $PA$ and $PB$ to the ellipse,   and let $x-y$ and $x+y$ be the coordinates of the points $A$ and $B$, where $0\leq y < 1/4$. Then $(x,y)$ are coordinates of the point $P$. We proved in Section \ref{back} that the coordinate curves $y=const$ and $x=const$ are  ellipses and hyperbolas, confocal with $\gamma$.

As in Section \ref{back}, the Poncelet grid is made by intersecting the tangent lines at points $\gamma(i/n), \ i=0,1,\dots,n-1$. The $(x,y)$-coordinates of the points of the grid are
$$
\left(\frac{k}{2n}+\frac{j}{n}, \frac{k}{2n}\right);\ k=0,1,\dots,\frac{n-1}{2},\ j=0,1,\dots,n-1.
$$
Fixing the second coordinate yields an angular set $P$ and fixing the first one -- a radial  set $Q$.

An ellipse 
$$
\frac{x_1^2}{a_1^2}+\frac{x_2^2}{a_2^2}=1
$$
also determines {\it elliptic coordinates} in the plane. Through a point $P$ there passes a unique ellipse and a unique hyperbola from the confocal family of conics (\ref{conffam}). The elliptic coordinates of $P$ are the respective values of the parameter, $\lambda_1$ and $\lambda_2$.  The  hyperbolas and ellipses from the confocal family (\ref{conffam}) are the coordinate curves of this coordinate system, $\lambda_1=const$ and $\lambda_2=const$, respectively. Cartesian coordinates of point $P$ are expressed in terms of the elliptic ones as follows:
\begin{equation} \label{Cartell}
x_1^2=\frac{(a_1^2+\lambda_1)(a_1^2+\lambda_2)}{a_1^2-a_2^2},\ x_2^2=\frac{(a_2^2+\lambda_1)(a_2^2+\lambda_2)}{a_2^2-a_1^2}
\end{equation}
(the Cartesian coordinates are determined up to the symmetries of an ellipse: $(x_1,x_2)\mapsto (\pm x_1,\pm x_2)$).

Thus the coordinates $(x,y)$ and the elliptic coordinates $(\lambda_1,\lambda_2)$ have the same coordinate curves, a family of confocal ellipses and hyperbolas. It follows that $\lambda_1$ is a function of $x$, and  $\lambda_2$ of $y$. 

Let $\Gamma_{\lambda}$ and $\Gamma_{\mu}$ be two ellipses (or two hyperbolas) from a confocal family of conics (\ref{conffam}). Consider the linear map 
$$
A_{\lambda,\mu} = {\rm Diag} \left(\sqrt{\frac{a_1^2 +\mu}{a_1^2 +\lambda}}, \sqrt{\frac{a_2^2 +\mu}{a_2^2 +\lambda}} \right).
$$
This map takes $\Gamma_{\lambda}$ to $\Gamma_{\mu}$. The following lemma is classical and goes back to J. Ivory.

\begin{lemma} \label{ivory}
If $\Gamma_{\lambda}$ and $\Gamma_{\mu}$ are two ellipses (respectively, two hyperbolas) and $P$ is a point of $\Gamma_{\lambda}$ then the points $P$ and $Q=A_{\lambda,\mu}(P)$ lie on the same confocal hyperbola (resp., ellipse).
\end{lemma} 

\noindent{\bf Proof}. We will argue in the case when $\Gamma_{\lambda}$ and $\Gamma_{\mu}$ are  ellipses. Let $(\lambda_1,\lambda_2)$ and $(\mu_1, \mu_2)$ be the elliptic coordinates of points $P$ and $Q$. Then $\lambda_2=\lambda$ and $\mu_2=\mu$. We want to prove that $\lambda_1=\mu_1$.

Let $(x_1,x_2)$ and $(X_1,X_2)$ be the Cartesian coordinates of points $P$ and $Q$.
One has formulas (\ref{Cartell}) and the similar relations:
\begin{equation} \label{Cartell2}
X_1^2=\frac{(a_1^2+\mu_1)(a_1^2+\mu_2)}{a_1^2-a_2^2},\ X_2^2=\frac{(a_2^2+\mu_1)(a_2^2+\mu_2)}{a_2^2-a_1^2}.
\end{equation}
On the other hand, $Q=A_{\lambda,\mu}(P)$, hence
$$
X_1^2=\frac{a_1^2 +\mu}{a_1^2 +\lambda} x_1^2 = \frac{(a_1^2+\lambda_1)(a_1^2+\mu_2)}{a_1^2-a_2^2},
$$
and likewise for $X_2^2$. Combined with (\ref{Cartell2}), this yields $\lambda_1=\mu_1$, as claimed. 
\proofend
  
Now we can prove that the sets $P_k$ and $P_m$
are linearly equivalent; the equivalence is given by the maps $\pm  A_{\lambda,\mu}$, depending on whether $k-m$ is even or odd. The argument for the sets $Q_k$ is similar.

The $(x,y)$-coordinates of the sets $P_k$ and $P_m$ are
$$
\left(\frac{k}{2n}+\frac{j}{n}, \frac{k}{2n}\right)\quad {\rm and}\quad 
\left(\frac{m}{2n}+\frac{j}{n}, \frac{m}{2n}\right);\ j=0,1,\dots,n-1.
$$
The sets $P_k$ and $P_m$ lie on confocal ellipses $\Gamma_{\lambda}$ and $\Gamma_{\mu}$. According to Lemma \ref{ivory}, the map  $A_{\lambda,\mu}$ preserves the first elliptic coordinate, and therefore, the $x$-coordinate. Therefore the coordinates of the points of the set $A_{\lambda,\mu} (P_k)$ are
$$
\left(\frac{k}{2n}+\frac{j}{n}, \frac{m}{2n}\right);\ j=0,1,\dots,n-1.
$$
If $m$ has the same parity as $k$, this coincides with the set $P_m$, and if the parity is opposite then this set is centrally symmetric to the set $P_m$. 

\section{Appendix 1: proof of integrability of the billiard in an ellipse} \label{int}
\begin{figure}[ht]
\centerline{\epsfbox{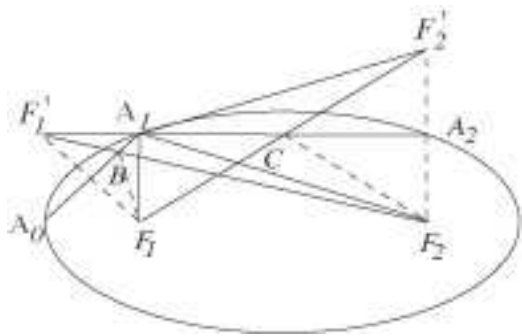}}
\caption{Integrability of the billiard in an ellipse}
\label{ellproof}
\end{figure}

Let $A_0 A_1$ and $A_1 A_2$ be consecutive segments of a billiard trajectory. Assume
that $A_0 A_1$ does not intersect the segment $F_1 F_2$; the other case is dealt with
similarly. It follows from the optical property of an ellipse, that the
angles $A_0 A_1 F_1$ and $A_2 A_1 F_2$ are equal; see figure \ref{ellproof}.

Reflect $F_1$ in $A_0 A_1$ to $F_{1}'$, and $F_2$ in $A_1 A_2$ to
$F_{2}'$, and set: $B = F_{1}' F_2 \cap A_0 A_1,  C = F_{2}' F_1 \cap A_1 A_2$.
Consider the ellipse with  foci $F_1$ and $F_2$ that is tangent to
 $A_0 A_1$. Since the angles $F_2 B A_1$ and $F_1 B A_0$ are equal, this ellipse
touches $A_0 A_1$ at the point $B$.  Likewise an ellipse with  foci $F_1$ and $F_2$
touches $A_1 A_2$ at the point $C$. One wants to  show that these two ellipses
coincide or, equivalently, that $F_1 B + B F_2 = F_1 C + C F_2$,  which boils down
to $F_{1}' F_2 = F_1 F_{2}'$.

Note that the triangles $F_{1}' A_1 F_2$ and $F_1 A_1 F_{2}'$
are congruent; indeed, $F_{1}' A_1 =  F_1 A_1,  F_2 A_1 = F_{2}' A_1$ by symmetry, and the angles $F_{1}' A_1 F_2$ and $F_1 A_1 F_{2}'$  are equal. Hence $F_{1}' F_2 = F_1 F_{2}'$, and the result follows.

\section{Appendix 2: the billiard map preserves measure $dS=\sin   \alpha \;  ds \wedge d \alpha $,  a geometrical explanation} \label{expl}

We provide a geometrical explanation of the area preserving property of the billiard ball map, Theorem \ref{area}.

\begin{figure}
\centerline{\epsfbox{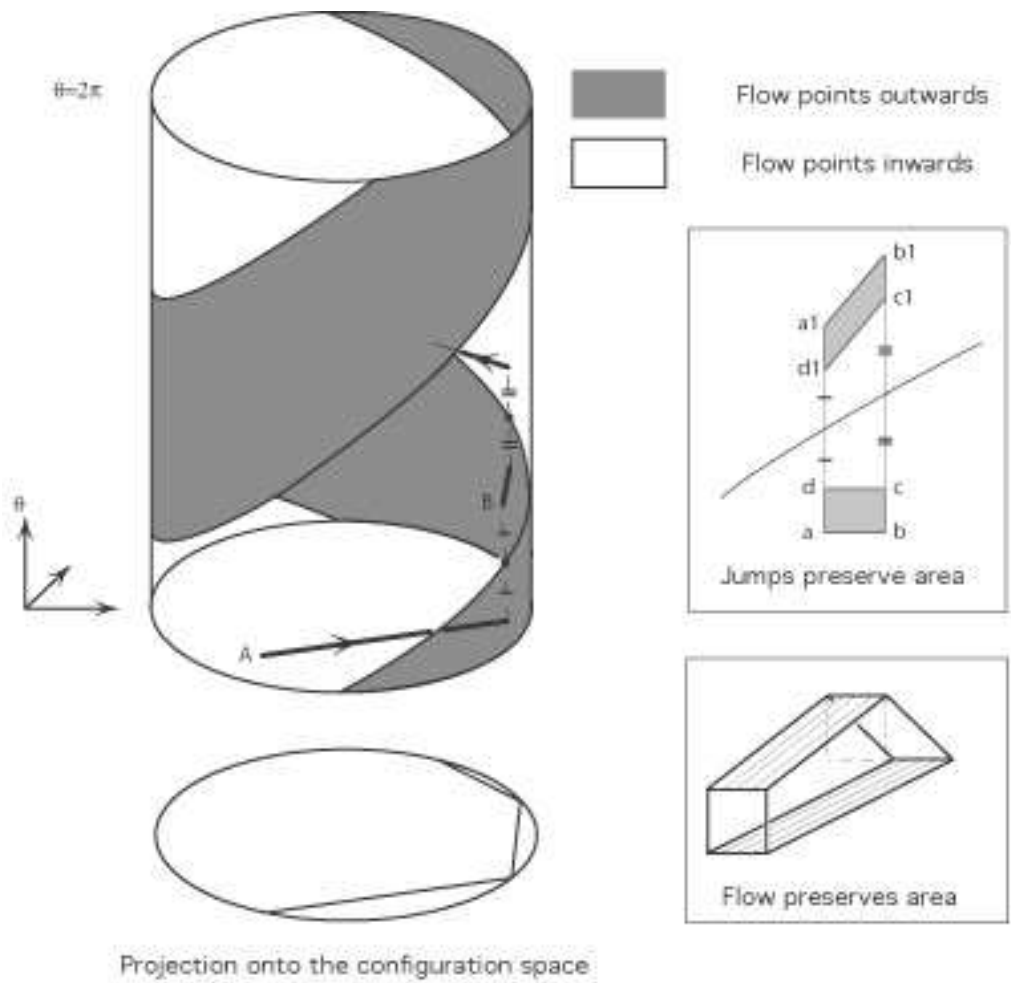}}
\caption{The billiard phase space}
\label{cyl}
\end{figure}

 {\it Properties of the billiard flow:} 
The phase space  of a billiard  (see Figure \ref{cyl}) consists of triples $( x, y, \theta ) $, where  $\theta$ is the angle between the particle's velocity and the $x$--axis. The flow is completely defined by the following two properties. 

\begin{enumerate} 
\item Each level $ \theta = \hbox{const.} $ carries a rigid translation in the $\theta$--direction\footnote{This reminds of airline traffic flow: collisions will be avoided if the planes at the same altitude   are at rest relative to each other.}: 
$ \dot x = \cos  \theta , \ \dot y  = \sin  \theta $. In particular, the flow inside the cylinder is volume-preserving (with the standard volume form  $ dx \wedge dy \wedge d\theta $).
\item  For any trajectory hitting  the boundary (the shaded part in the Figure \ref{cyl}),  $\theta$ jumps according to  
the incidence-reflection law $ \theta ^\prime- \beta   = - (\theta - \beta ) $,  where $\beta  = \beta (s)  $ is the angle of the tangent to the boundary at the collision point with the $x$--axis, and $s$ is the length parameter along the boundary of the billiard table. This jump preserves the area form $ ds\wedge d\theta$ of the boundary to itself -- indeed, it is a reflection in $\theta$ around $ \beta = \beta (s) $ for each $s$. 
\end{enumerate} 

\begin{figure}[ht]
\centerline{\epsfbox{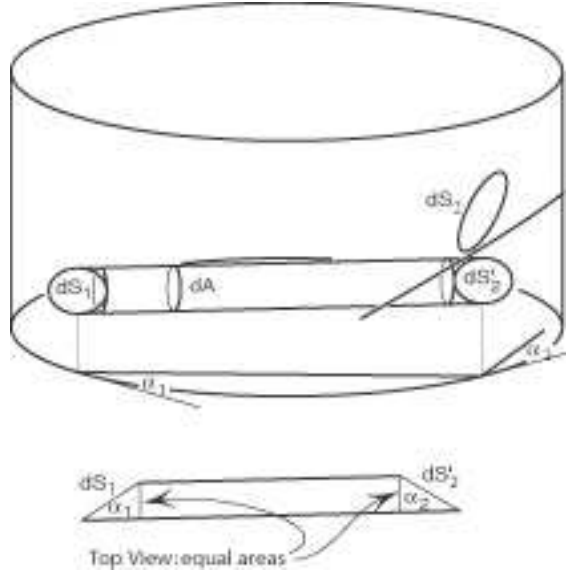}}
\caption{$ dS_1\sin \alpha _1 = dS_2 ^\prime  \sin \alpha _2 =dS_2  \sin \alpha _2  $}
\label{last}
\end{figure}

 {\it Explanation of the area--preservation:}
Consider the Poincar\'e map which takes points $A$ to $B$ as shown in figure. 
The map is defined on the (unshaded) region where the flow is directed inwards: 
$ \beta (s) \leq \theta \leq \beta (s) + \pi $. This map preserves  
$\sin \alpha\ ds\wedge d\theta$ , 
as explained in  Figure \ref{last}. Indeed, consider the flow tube originating on a small patch around $A$. The normal cross-sectional area $ dA$   is constant along the tube, since  the flow is volume-preserving and has unit speed. But  for the areas of the oblique cross-sections we have  $ dS_1= dA/ \sin \alpha _1, \ dS_2 ^\prime  = dA/ \sin \alpha _2 $. Since $ dS_2^\prime =dS_2 $ as shown before, we have 
$\sin \alpha\ dS_1 = \sin \alpha\ dS_2$
It remains to observe that the area form on the boundary is $ ds\wedge d\theta = 
ds\wedge d\alpha $. 
 \bigskip

{\bf Acknowledgments.} Many thanks to R. Schwartz for  fruitful discussions of his beautiful theorem. We are grateful to the Mathematics Institute at Oberwolfach for its hospitality. Both authors were supported in part by NSF.

\end{document}